\documentclass[12pt]{amsart}

\usepackage{amsmath,amssymb,amsthm,amscd, eucal}
\usepackage{amscd,enumerate} 
\usepackage{mathtools} 
\usepackage{latexsym}
\usepackage{bm} 
\usepackage[usenames]{color}
\usepackage[utf8]{inputenc} 
\usepackage{dsfont}
\usepackage{graphicx} 

\usepackage[shortlabels]{enumitem} 
\setlist{nolistsep} 

\newtheorem{thm}[equation]{Theorem}

\newtheorem{cor}[equation]{Corollary}
\newtheorem{lem}[equation]{Lemma}

\newtheorem{prop}[equation]{Proposition}

\theoremstyle{definition}

\newtheorem{defn}[equation]{Definition}

\newtheorem{defns}[equation]{Definitions}
\newtheorem{exer}{Exercise}
\newtheorem*{sol}{Solution}
\newtheorem{exam}[equation]{Example}

\newtheorem{rem}[equation]{Remark}

\numberwithin{equation}{section}

\newcommand{\eq}[1]{\begin{align*}#1\end{align*}} 

\newcommand{\id}{{\rm id}} 

\newcommand{\ria}{\rightarrow}

\newcommand{\n}{\textbf}
\newcommand{\e}{\emph}

\newcommand{\enum}{\begin{enumerate}[leftmargin=*,itemsep=-0.5ex]}
\newcommand{\enumi}{\begin{enumerate}[{\rm i)},leftmargin=*,itemsep=-0.5ex]}
\newcommand{\enuma}{\begin{enumerate}[{\rm a)},leftmargin=*,itemsep=-0.5ex]}
\newcommand{\eenum}{\end{enumerate}}
\newcommand{\ite}{\begin{itemize}[leftmargin=*,itemsep=-0.5ex]}
\newcommand{\eite}{\end{itemize}}

\newcommand{\ce}{\begin{center}}
\newcommand{\ece}{\end{center}}
\newcommand{\bmp}{\begin{minipage}}
\newcommand{\ebmp}{\end{minipage}}
\newcommand{\bma}{\begin{matrix}}
\newcommand{\ema}{\end{matrix}}

\newcommand{\de}{\begin{defn}}
\newcommand{\ede}{\end{defn}}
\newcommand{\pr}{\begin{prop}}
\newcommand{\epr}{\end{prop}}
\newcommand{\tm}{\begin{thm}}
\newcommand{\etm}{\end{thm}}
\newcommand{\lm}{\begin{lem}}
\newcommand{\elm}{\end{lem}}
\newcommand{\co}{\begin{cor}}
\newcommand{\eco}{\end{cor}}

\newcommand{\p}{\begin{proof}}
\newcommand{\ep}{\end{proof}}
\newcommand{\ej}{\begin{exer}}
\newcommand{\eej}{\end{exer}}
\newcommand{\so}{\begin{sol}}
\newcommand{\eso}{\end{sol}}
\newcommand{\ex}{\begin{exam}}
\newcommand{\eex}{\end{exam}}
\newcommand{\re}{\begin{rem}}
\newcommand{\ere}{\end{rem}}

\newcommand{\bi}{\begin{thebibliography}{}}
\newcommand{\ebi}{\end{thebibliography}}

\begin{document}

\title[M(G)S derivations in prime rings]{There are no nontrivial two-sided multiplicative (generalized)-skew derivations in prime rings}

\author{Jose Brox}
\address{University of Coimbra, CMUC, Department of Mathematics, 3004-504 Coimbra, Portugal}
\thanks{This work was partially supported by the Centre for Mathematics of the University of Coimbra - UIDB/00324/2020, funded by the Portuguese Government through FCT/MCTES. The author was supported by the Portuguese Government through the FCT grant SFRH/BPD/118665/2016.}
\email{josebrox@mat.uc.pt}


\keywords{Prime rings, derivations, generalized derivations, skew derivations}
\subjclass[2010]{16N60, 16W25, 16W20}

\maketitle

\begin{abstract}
As originally defined by Ashraf and Mozumder (\cite{AM}), multiplicative (generalized)-skew derivations must satisfy two identities. In this short note we show that, as a consequence of the simultaneous satisfaction of both identities, a multiplicative (generalized)-skew derivation of a prime ring is either a multiplicative (generalized) derivation (i.e., not skew), or a generalized skew derivation (i.e., additive). Therefore only one of the identities should be taken in the definition of multiplicative (generalized)-skew derivations in order to get a new class of derivations in prime rings.
\end{abstract}

\section{Introduction}

The fundamental concept of \e{derivation} of an associative ring $R$, an additive map $d:R\ria R$ such that $d(xy)=d(x)y+xd(y)$, has been progressively generalized in recent literature: by a twisting by an automorphism or a secondary derivation of the ring, by dropping the additivity assumption, by combining both previous ideas, and by repeating the process on the secondary derivation when present (\cite{D,B,L,DT,DA,AM}). One of the main purposes of these generalizations is to extend to more sophisticated maps the classic results on derivations in the tradition of Herstein's theory of rings (\cite{MM}), in which strong knowledge is gained about the map or the ring through some special (and a priori weaker) property of the map. The main focus is on prime and semiprime rings, or on rings with well-behaved idempotents, which provide a context rich enough for the theory to be satisfactorily developed. For example one tries to extend to a more general setting Posner's second theorem for derivations of prime rings (\cite[Theorem 2]{P}), which states that a prime ring $R$ is commutative when it has a derivation $d\neq0$ such that $xd(x)-d(x)x$ is central for every $x\in R$. These efforts have generated literature in abundance (e.g. \cite{AMB,C,FRR,DP,Hu,H,KG,LS,RK}).

\begin{defns}
Let $R$ be a ring.

\enuma

\item A \e{skew derivation} (\cite[page 170]{J}) is an additive map $d:R\ria R$ together with an automorphism $\alpha:R\ria R$ such that $d(xy)=d(x)y+\alpha(x)d(y)$.
\item A \e{multiplicative derivation} (\cite{D}) is a map $d:R\ria R$, not necessarily additive, such that $d(xy)=d(x)y+xd(y)$.
\item We define a \e{multiplicative skew derivation} as a not necessarily additive map $d:R\ria R$ together with an automorphism $\alpha:R\ria R$ such that $d(xy)=d(x)y+\alpha(x)d(y)$.
\item A \e{generalized derivation} (\cite{B}) is an additive map $F:R\ria R$ together with a derivation $d:R\ria R$ such that $F(xy)=F(x)y+xd(y)$.
\item A \e{generalized skew derivation} (\cite{L}) is an additive map $F:R\ria R$ together with an automorphism $\alpha:R\ria R$ and a skew derivation $d:R\ria R$ for $\alpha$ such that $F(xy)=F(x)y+\alpha(x)d(y)$.
\item  A \e{multiplicative generalized derivation} (\cite{DT}) is a map $F:R\ria R$, not necessarily additive, together with a derivation $d:R\ria R$ such that $F(xy)=F(x)y+xd(y)$.
\item  A \e{multiplicative (generalized) derivation} (\cite{DA}) is a map $F:R\ria R$, not necessarily additive, together with a map (not necessarily additive nor a derivation) $d:R\ria R$ such that $F(xy)=F(x)y+xd(y)$.
\eenum
As defined in \cite{AM} (see also \cite{MD}), a \e{multiplicative (generalized)-skew derivation} (M(G)S derivation) is a not necessarily additive map $F:R\ria R$, together with a not necessarily additive map $d:R\ria R$ and an automorphism $\alpha:R\ria R$, satisfying the identities
\eq{F(xy)=&F(x)y+\alpha(x)d(y) \ \text{ \n{(Identity 1)}}\\
F(xy)=&F(x)\alpha(y)+xd(y) \ \text{ \n{(Identity 2)}}}
Since in this case we get two different identities in the definition, accordingly we will call these maps \e{two-sided M(G)S derivations}. We will say that a map is a \e{M(G)S derivation of type 1} (resp. \e{M(G)S derivation of type 2}) when it satisfies Identity 1 (resp. Identity 2).
\end{defns}

\section{Main theorem}

In what follows we show that in prime rings there are no nontrivial two-sided M(G)S derivations, since each of them is either a multiplicative (generalized) derivation (i.e., not skew) or a generalized skew derivation (i.e., additive with additive derivation). Therefore, in prime rings the only interesting M(G)S derivations are either of type 1 or of type 2.

\lm\label{lemma}
If $R$ is a semiprime ring and $F$ is a M(G)S derivation of type 1 (resp. type 2) with map $d:R\ria R$ and automorphism $\alpha:R\ria R$ then $d$ is in fact a multiplicative skew derivation with $\alpha$ as automorphism (resp. satisfies $d(xy)=d(x)\alpha(y)+xd(y))$.
\elm

\p
For type 1 this is \cite[Lemma 2.1]{RK}, reproduced here for convenience. If $F(xy)=F(x)y+\alpha(x)d(y)$ (1) then on the one hand \[F(xyz)=F(x)yz+\alpha(x)d(yz)\]
for all $x,y,z\in R$, while on the other
\[F(xyz)=F(xy)z+\alpha(xy)d(z)\overset{(1)}{=}F(x)yz+\alpha(x)d(y)z+\alpha(x)\alpha(y)d(z),\]
hence $\alpha(x)d(yz)=\alpha(x)d(y)z+\alpha(x)\alpha(y)d(z)$ for all $x,y,z\in R$, and
\[\alpha(x)(d(yz)-d(y)z-\alpha(y)d(z))=0.\]
Since $\alpha$ is surjective we get $R(d(yz)-d(y)z-\alpha(y)d(z))=0$, so by semiprimeness $d(yz)=d(y)z+\alpha(y)d(z)$ for all $y,z\in R$.
For type 2 the proof is analogous.
\ep

\tm
Let $R$ be a prime ring and $F$ be a two-sided M(G)S derivation with map $d:R\ria R$ and automorphism $\alpha:R\ria R$. Then
\enumi
\item $\alpha=\id_R$, so $F$ is a multiplicative (generalized) derivation, or
\item $F$ and $d$ are additive, so $F$ is a generalized skew derivation.
\eenum
\etm

\p
From Identities 1 and 2, $F(x)y+\alpha(x)d(y)=F(xy)=F(x)\alpha(y)+xd(y)$ for every $x,y\in R$,
so \[F(x)(y-\alpha(y))=(x-\alpha(x))d(y).\]
Since the expression $(x-\alpha(x))d(y)$ is linear in $x$ for each $y$, so is $F(x)(y-\alpha(y))$. Hence
\[F(x+y)(z-\alpha(z)) = F(x)(z-\alpha(z))+F(y)(z-\alpha(z)).\]
So $(F(x+y)-F(x)-F(y))(z-\alpha(z))=0$ for every $x,y,z\in R$. Put $G(x,y):=F(x+y)-F(x)-F(y)$. We have, for every $x,y,z\in R$,
\[G(x,y)z=G(x,y)\alpha(z). \eqno{\n{(1)}}\]
Therefore, for every $w\in R$,
\[G(x,y)wz\overset{\n{(1)}}{=}G(x,y)\alpha(wz)=(G(x,y)\alpha(w))\alpha(z)\overset{\n{(1)}}{=}(G(x,y)w)\alpha(z),\]
hence
$G(x,y)w(z-\alpha(z))=0$ for every $x,y,z,w\in R$. Since $R$ is prime, either $\alpha(z)=z$ for every $z\in R$ or $G(x,y)=0$ for every $x,y\in R$. In the first case $\alpha=\id_R$ and $F$ is a multiplicative (generalized) derivation. In the second case we get $\alpha\neq\id_R$ and $F(x+y)=F(x)+F(y)$ for every $x,y\in R$, so $F$ is additive. Now, by Lemma \ref{lemma} above $d$ is another M(G)S derivation associated to $\alpha\neq\id_R$, so analogously $d$ is additive, whence it is a skew derivation and $F$ is a generalized skew derivation.
\ep

\bi
\bibitem{AMB} Ali, A., Muthana N., Bano A. (2018). Multiplicative (generalized)-derivations and left multipliers in semiprime rings.  \e{Palest. J. Math.} 7, Special Issue I:170--178.
\bibitem{AM} Ashraf, M., Mozumder, M.R. (2015). On multiplicative (generalized)-skew derivations over semiprime rings.  \e{Rend. Sem. Mat. Univ. Pol. Torino} 73/2(3-4):261--268.
\bibitem{B} Bre\u sar, M. (1991). On the distance of the composition of two derivations to the generalized derivations. \e{Glasgow Math. J.} 33:89--93.
\bibitem{C} Chang, J.-C. (2010). Generalized skew derivations with nilpotent values on Lie ideals. \e{Monatsh. Math.} 161(2):155--160.
\bibitem{D}  Daif, M.N. (1991). When is a multiplicative derivation additive? \e{Int. J. Math. Math. Sci.} 14(3):615--618.
\bibitem{DT} Daif, M.N., Tammam El-Sayiad, M.S. (2007). Multiplicative generalized derivations which are additive. \e{East-West J. Math.} 9(1):31--37.
\bibitem{FRR} De Filippis, V., Rehman, N., Raza, M.A. (2018). Strong commutativity preserving skew derivations in semiprime rings.  \e{Bull. Malays. Math. Sci. Soc.} 41(4):1819--1834.
\bibitem{DA} Dhara, B., Ali, S. (2013). On multiplicative (generalized)-derivations in prime and semiprime rings. \e{Aequationes Math.}
86(12):65–-79.
\bibitem{DP} Dhara, B., Pradhan, K.G. (2016). A note on multiplicative (generalized) derivations with annihilator conditions. \e{Georgian Math. J.} 23(2):191--198.
\bibitem{Hu} Huang, S. (2018). Commuting multiplicative generalized derivations on Lie ideals of semiprime rings. \e{Gulf J. Math.} 6(3):25--32.
\bibitem{H} Hvala, B. (1998). Generalized derivations in rings. \e{Comm. Algebra} 26(4):1147--1166.
\bibitem{J} Jacobson, N. (1964). \e{Structure of Rings}. Providence: American Mathematical Society. Revised edition.
\bibitem{KG} Koç, E., Gölba\c{s}i, Ö. (2018). Some results on ideals of semiprime rings with multiplicative generalized derivations.  \e{Comm. Algebra} 46(11):4905--4913.
\bibitem{L} Lee, T.K. (2004). Generalized skew derivations characterized by acting on zero products. \e{Pacific J. Math.} 216(2):293--301.
\bibitem{LS} Lee, T.K., Shiue, W.-K. (2001). Identities with generalized derivations. \e{Comm. Algebra} 29(10):4437--4450.
\bibitem{MD} Mozumder, M.R., Dhara, B. (2016). Multiplicative (generalized)-skew derivations on semiprime rings. In:\e{International conference on semigroups, algebras and applications.}
\bibitem{MM} Martindale III, W.S., Miers, C.R. (1986). Herstein's Lie theory revisited. \e{J. Algebra} 98(1):14--37.
\bibitem{P} Posner, E.C. (1957). Derivations in prime rings. \e{Proc. Amer. Math. Soc.} 8:1093--1100.
\bibitem{RK} Rehman, N., Khan, M.S. (2018). A note on multiplicative (generalized)-skew derivation on semiprime rings. \e{Journal of Taibah University for Science} 12(4):450--454.
\ebi

\end{document}